\documentclass[pre,preprint,showpacs,amsmath,amssymb]{revtex4-1}

\usepackage{mathtools}
\usepackage{amsmath}
\usepackage{amsfonts}
\usepackage{amssymb}
\usepackage{graphicx}
\usepackage{epstopdf}
\usepackage{dsfont}
\usepackage{amsthm}
\usepackage{hyperref}
\hypersetup{
    colorlinks = false
}
\usepackage{float}
\usepackage{xfrac}
\usepackage{cancel}

\newcommand*{\eqb}{\begin{eqnarray}}
\newcommand*{\eqe}{\end{eqnarray}}
\newcommand*{\al}{\alpha}

\DeclareMathOperator{\D}{\operatorname{d}\!}
\DeclareMathOperator{\E}{\mathbb{E}}
\DeclareMathOperator{\PR}{\mathbb{P}}
\DeclareMathOperator{\R}{\mathbb{R}}
\DeclareMathOperator{\F}{\mathcal{F}}

\newcommand{\sgn}{\operatorname{sgn}}

\theoremstyle{plain}

\numberwithin{theorem}{section}

\theoremstyle{definition}
\numberwithin{theorem}{section}

\begin{document}

\title{Stochastic representation of processes with resetting}

\author{Marcin Magdziarz}

\email[corresponding author: ]{marcin.magdziarz@pwr.edu.pl}

\affiliation{Faculty of Pure and Applied Mathematics, Hugo Steinhaus Center, Wroclaw University of Science and Technology, Wyspianskiego 27,
50-370 Wroclaw, Poland}
\author{Kacper Ta\'zbierski}

\affiliation{Faculty of Pure and Applied Mathematics, Hugo Steinhaus Center, Wroclaw University of Science and Technology, Wyspianskiego 27,
50-370 Wroclaw, Poland}

\date{September 7, 2021}

\begin{abstract}
In this paper we introduce a general stochastic representation for an important class of processes with resetting. It allows to describe any stochastic process intermittently terminated and restarted from a predefined random or non-random point. Our approach is based on stochastic differential equations called jump-diffusion models. It allows to analyze processes with resetting both, analytically and using Monte Carlo simulation methods. To depict the strength of our approach, we derive a number of fundamental properties  of Brownian motion with Poissonian resetting, such as: the It\^o lemma, the moment-generating function, the characteristic function, the explicit form of the probability density function, moments of all orders, various forms of the Fokker-Planck equation, infinitesimal generator of the process and its adjoint operator. Additionally, we extend the above results to the case of time-nonhomogeneous Poissonian resetting.  This way we build a general framework for the analysis of any stochastic process with intermittent random resetting.
\end{abstract}

\maketitle

\section{Introduction}
In recent years we have been observing a great interest in the field of intermittent stochastic processes, in which the
motion of a particle is interrupted by random resetting to an initial state.
Random motion with resetting is typically observed in various searching strategies and foraging patterns, where 
the living organism after an unsuccessful excursion returns to the initial position and start the search again \cite{campos2014reorientation}.  
The justification of such strategy is that if one does
not succeed in finding the target within certain not very long time period, sometimes it
is more beneficial and safe to move back to the origin and start another excursion, see 
\cite{manrubia1999stochastic,gelenbe2010search}.
On the other hand, in the field of population dynamics resetting could be interpreted as an event reducing population to its natural size according to the environmental capacity \cite{Kyriakidis1994}. What makes the idea also interesting is the way that random resetting completely changes the properties of the diffusion process. In particular,  mean first passage time of diffusion with resetting is finite, which is in sharp contrast to the standard diffusion case \citep{evans2011diffusion,Evans2020}. This fact depicts the efficiency of searching strategies with resetting mechanism \cite{Tong2008}. 
However, if the resetting timer has an infinite mean then mean first passage time may be infinite.
Intermittent processes have found important application in other fields, such as: 
phenotypic diversity, population growth and information in fluctuating environments \cite{kussell2005phenotypic}, ecology \cite{viswanathan2011physics}, enzymatic reactions \cite{reuveni2014role} or computer science and optimization problems  \cite{montanari2002optimizing}.
The historical background of processes with resetting can be found in \cite{montero2017continuous}. The interested reader is referred to a recent comprehensive review by Evans et al. \cite{Evans2020}.
Some latest results related to non-Markovian resetting models can be found in \cite{bodrova2019nonrenewal,zhou2020continuous,zhou2021gaussian}. 
We also refer the reader to the following relevant works related to the restarted processes, published well before the terms "restart" and "resetting" became popular \cite{Iddo1,Iddo2,Iddo3}.

To characterize a process with resetting, one needs to deal with two sources of randomness. The first one accounts for the particle motion between resetting times, the second one is the point process determining the moments of resetting.
The usual description of motion with resetting is via CTRW with additional resetting mechanism or via related deterministic partial differential equations of the Fokker-Planck type \citep{Evans2020}. 

Here we use another approach, which is founded on the theory of L\'evy processes and the corresponding stochastic differential equations (SDEs) -- the so-called jump-diffusion models.
There is a rich history of research of  L\'evy processes. The most well-known is obviously Brownian motion  (or Wiener process). Other important examples include $\al$-stable, Linnik, Mittag-Leffler, Gamma, or Laplace processes \cite{ken1999levy}.
Their first applications arose quickly in the physical \cite{Einstein1905,Smoluchowski1906} and financial \cite{Bachelier,merton1976,Black-Scholes} sciences. They have naturally spread to many other fields of science, e.g. as models in biology, chemistry, data mining, statistics, etc. \cite{barndorff2012levy}.
In 1940's, thanks to Kiyoshi It\^o, the theory of stochastic differential equations was developed \cite{Ito1}, \cite{Ito2}. The stochastic It\^o integral and the celebrated It\^o lemma allowed to integrate functions or processes with respect to the Wiener process. Also, other types of stochastic integrals with important role in physics, such as Stratanovich integral, were developed. 
Jump-diffusion models are defined as SDEs, in which the Brownian noise term is replaced or complemented by arbitrary L\'evy noise \cite{OksendalSulem}. Since general L\'evy processes have discontinuous trajectories \cite{ken1999levy}, therefore solutions of such L\'evy driven SDEs display jumps.
Jump-diffusion models have found applications in various fields, such as: finance, statistical physics, pattern theory and computational vision and many more, see \cite{OksendalSulem,eliazar2003levy} and references therein. As we will show in the following, these processes are also tailor-made to describe and analyze processes with resetting.

In this paper we introduce a general stochastic framework for processes with resetting in the language of stochastic differential equations (SDEs) with jumps (jump-diffusion models). Using the stochastic approach a version of the It\^o lemma for processes with resetting will be derived.  Next, applying the latter result we will  derive the moment generating function (MGF) and the Fourier transform for diffusion with resetting.  Using the inverse Fourier transform we will then determine the probability density function (PDF) of the process in explicit form. The corresponding Fokker-Planck equation will be also derived in various forms. We will discuss the equivalence of our Fokker-Planck equation with the one given in \cite{Evans2020}. Moreover we will find the infinitesimal generator and the adjoint operator for the Markovian resetting process. In addition, we will extend the above results to the case of time-nonhomogeneous Poissonian resetting. We will show that, depending on the intensity function of the resetting events, we can observe: convergence to stationary distribution, convergence to a point, subdiffusive or diffusive behaviour, for long times.
 
It should be underlined that the presented here approach can be easily generalized to the case of non-Markovian resetting times as well as other, not necessarily Brownian, driving processes (L\'evy processes, fractional processes, etc.).  

\section{Stochastic representation of processes with resetting}

\begin{figure}[hb]
    \includegraphics[scale=0.3]{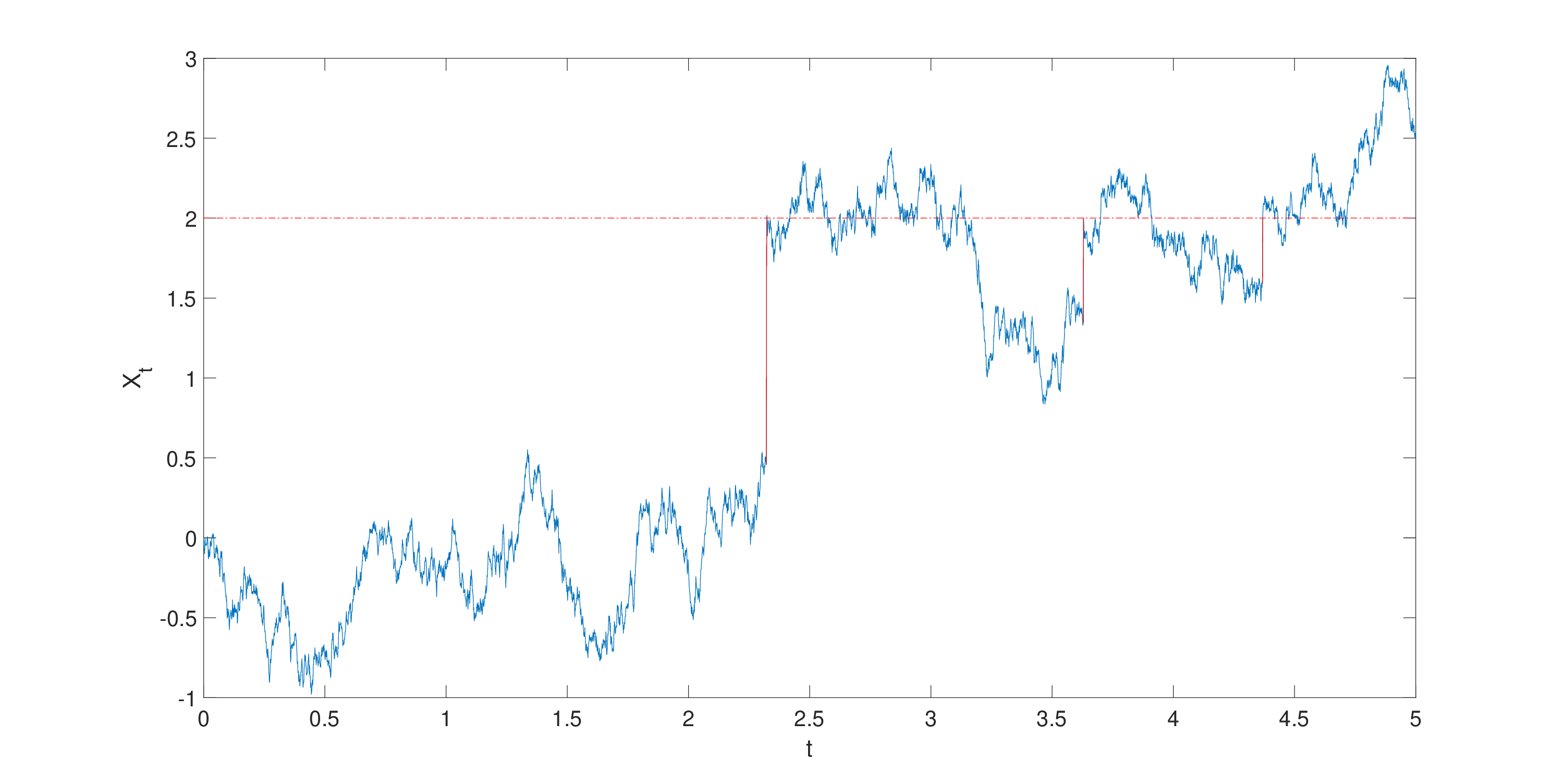}
    \caption{An exemplary plot of a resetting process. The process starts at $x_0=0$ and then diffuses normally until it moves instantaneously to $x_{R}=2$ in every consecutive resetting event.}
    \label{fig:resettingplot}
\end{figure}

Let us consider a diffusing particle with initial position $x_0$ and Poissonian resetting with rate $r$
to position $x_{R}$. In general the initial position $x_0$ and the resetting position $x_{R}$ can be distinct, they can be even random  (in this paper we assume that $x_0$ and $x_R$ are constant, however the methodology we use, can also be applied to the random case, which will lead to different results and distributions). Position of such particle at time $t$ is usually defined in the following way \cite{Evans2020}:
\begin{equation}\label{eq:evansSDE}
\begin{gathered}
    x(t+\D t)=x_{R}\quad \text{with probability }r\D t\\
    =x(t)+\xi (t)(\D t)^{\frac{1}{2}}\quad \text{with probability }(1-r\D t),
    \end{gathered}
\end{equation}
where $\xi(t)$ is the standard Gaussian white noise with mean 0 and variance $2D$. 
The typical trajectory of a diffusion process with resetting is depicted in Fig. \ref{fig:resettingplot}.
The Fokker-Planck equation corresponding to \eqref{eq:evansSDE} has the form \cite{Evans2020}
\begin{equation}\label{eq:evanspdf}
    \partial_tp(x,t)=D \partial_{xx}p(x,t)+r\delta(x-x_{R})-rp(x,t),
\end{equation}
with initial condition $p(x,0) = \delta(x - x_0)$. Here $p(x,t)$ is the PDF of the diffusing particle with resetting and $\delta(\cdot)$ is the Dirac delta. In what follows we will put for simplicity $D=1/2$. Using the above formalism one can examine the interesting phenomenon of non-equilibrium steady states \cite{Evans2020}. Even though the stationary state is independent of time, there is a driving force in form of the resetting that creates the probability flow. In a purely mathematical setting it means that the gradient of stationary state is non-zero. This corresponds to the physical concept of non-equilibrium steady state. A nice review on this topic can be found in \cite{ge2012stochastic}.

%
 
Also in \cite{Evans2020} an absorption of the process by a trap was studied. For the standard diffusion, time to absorption  follows L\'evy distribution, which has an infinite mean (being an $\alpha$-stable distribution with $\alpha=1/2$). When the resetting is introduced, the mean time to absorption becomes finite and we are able to find an optimal resetting rate \cite{Evans2020}. This fact clearly shows that resetting can be beneficial in the searching strategies.


In this paper we introduce a different approach to define and analyze processes with resetting. Namely, we will use jump-diffusion processes to define the stochastic dynamics with resetting.
A jump-diffusion process is defined as the solution of the following SDE \cite{OksendalSulem}:
\eqb
\label{jump_diffusion}
dX_t = \mu_t dt +\sigma_t dW_t + \nu_t dL_t \; , \;\;\; X_0 = x_0.
\eqe
Here $W_t$ is the standard Brownian motion, $L_t$ the L\'evy process, it is usually the Poisson process and introduces jumps to the observed dynamics, $\mu_t$, $\sigma_t$ and $\nu_t$ are the appropriate parameters of the model, which can be in general space and time dependent. It is also assumed that $W_t$ and $L_t$ are independent.
Recall that L\'evy process $L_t$ is a stochastic process satisfying \cite{OksendalSulem}:
\begin{enumerate}
    \item $L_0=0$,
    \item  $L_t$ has independent increments,
    \item $L_t$ has stationary increments,
    \item $\displaystyle\mathop{\forall}_{\varepsilon,t>0}\lim_{h\to 0}P\left( |L_{t+h}-L_t|>\varepsilon \right)=0$, i.e. $L_t$ is continuous in probability.
\end{enumerate}
Since every L\'evy process is a Markov process, solutions of SDE \eqref{jump_diffusion}  are Markovian.
Typical examples of L\'evy processes are: Brownian motion, Poisson process,  and $\alpha$-stable L\'evy process.
Brownian motion is the only L\'evy process (up to a constant drift) with continuous trajectories. All the other L\'evy processes have jumps. 
Poisson process is a non-decreasing L\'evy process with jumps of size $1$ and flat periods between jumps. Times between consecutive jumps of the Poisson process are independent and drawn from the exponential distribution with mean $r>0$. The constant $r$  is called intensity of the Poisson process

Now let us consider the following particular case of jump-diffusion process:
\begin{equation}\label{eq:mysde}
    \D X_t=\D W_t+\left( x_{R}-X_t \right)\D N_t,\quad X_0=x_0,
\end{equation}
where $N_t$ is the Poisson Process with intensity $r$, $x_0$ and  $x_{R}$ are constants. 
We argue that \eqref{eq:mysde} is the continuous-time analogue of diffusion with Poissonian resetting defined in \eqref{eq:evansSDE}.
Indeed, whenever $N_t$ has jump, we get that $dN_t =1$ and therefore the term $\left( x_{R}-X_t \right)\D N_t$ 
sends back the process to the resetting position $x_{R}$.
Between the resetting times, i.e. when $N_t$ is constant, we get that $dN_t$=0. So in this case $X_t =W_t$. Thus between consecutive resetting times the particle performs Brownian motion. 
A typical trajectory obtained using \eqref{eq:mysde}  can be seen in Fig, \ref{fig:resettingplot}.
Clearly, we observe a diffusing particle with resetting to $x_{R} = 2$.

We can also write down the integrated form of \eqref{eq:mysde}:
\begin{equation}\label{eq:sdesol}
    X_t = x_0+W_t+ \int_0^t (x_{R}-X_s)dN_s=x_0+W_t+\sum\displaylimits_{n=1}^{N_t}\left( x_{R}-X_{\tau_n}\right).
\end{equation}
Here $\tau_n=\max \{ t>0: N_t\leq n\}$ are the consecutive resetting times. 

In what follows we will show that the derived here stochastic representation  \eqref{eq:mysde} constitutes a general framework for the analysis of processes with resetting. It allows to study many key properties of intermittent processes both analytically and using Monte Carlo simulation methods.
Moreover, it should be underlined that the introduced approach can be easily generalized to the case of non-Poissonian resetting times as well as other arbitrary driving processes. More precisely, suppose that we want to define a stochastic process 
$\tilde{X}_t$ describing an arbitrary dynamics with resetting. Let $T_1,T_2,...$ be an arbitrary sequence of positive random variables describing times between consecutive resetting events. Moreover, let us assume that $\tilde{W}_t$ is an arbitrary stochastic process describing the particle motion between resetting events. Then $\tilde{X}_t$ can be defined as the solution of the following SDE:
\begin{equation}\label{eq:mysde2}
    \D \tilde{X}_t=\D \tilde{W}_t+\left( x_{R}-\tilde{X}_t \right)\D \tilde{N}_t,\quad \tilde{X}_0={x}_0,
\end{equation}
which is a straightforward generalization of \eqref{eq:mysde}. Here
$\tilde{N}_t=\max\{n\in\mathbb{N}: \sum_{i=1}^n T_i \leq t\}$. It counts the number of resetting events up to time $t$. In particular, $\tilde{N}_t$ can be a renewal processes.
Consequently, any process with resetting can be written in the form of \eqref{eq:mysde2}.

\subsection{It\^o lemma}
It\^o lemma is the main tool in the analysis of SDEs. It allows to find their solutions as well as many of their key properties. 
Let us derive a version of the It\^o lemma corresponding to the resetting process \eqref{eq:mysde}. 
For completeness we will derive it in a more general setting.
Let $Z_t=f(X_t)$, where $f$ is appropriately smooth and $X_t$ is governed by the following SDE
\eqb
\label{jump_diffusion3}
dX_t = \mu_t dt +\sigma_t dW_t + (x_{R} - X_t)dN_t \; , \;\;\; X_0 = x_0.
\eqe
Note that for $\mu_t = 0$ and $\sigma_t =1$ we recover the resetting process \eqref{eq:mysde}.
Using the Taylor expansion we get
\begin{align*}
    \D f(X_t)=\sum\displaylimits_{i=1}^\infty\frac{f^{(i)}(X_t)}{i!}\left( \D X_t \right)^i.
\end{align*}
According to classical It\^o calculus we have that $(\D W_t)^2 = \D t$ and $(\D W_t)^n = 0$ for $n>2$. Moreover
$(\D N_t)^n = \D N_t$ for any $n\in\mathbb{N}$ since $N_t$ is the Poisson process. Therefore 
\begin{equation*}
   (\D X_t)^2 =\left( \sigma^2_t\D t+\left( x_{R}-X_t \right)^2\D N_t\right)
\end{equation*}
and
\begin{equation*}
(\D X_t)^i =\left( x_{R}-X_t \right)^i\D N_t
\end{equation*}
for $i>2$.
Consequently
\begin{align*}
    \D f(X_t)&= f^{'}(X_t)\left(\mu_t\D t + \sigma_t \D W_t+\left( x_{R}-X_t \right)\D N_t \right)+\\
    &+\frac{1}{2}f^{''}(X_t)\left( \sigma^2_t\D t+\left( x_{R}-X_t \right)^2\D N_t)\right)+\sum\displaylimits_{i=3}^\infty\frac{f^{(i)}(X_t)}{i!}\left(x_{R}-X_t\right)^i \D N_t \\
    &=f ^{'}(X_t)\sigma_t\D W_t+\left(f^{'}(X_t) \mu_t+f^{''}(X_t) \frac{\sigma_t^2}{2} \right)\D t+\\
    &+\left( f ^{'}(X_t)\left( x_{R}-X_t \right)+ \frac{1}{2}f^{''}(X_t)\left( x_{R}-X_t \right)^2+\sum\displaylimits_{i=3}^\infty\frac{f^{(i)}(X_t)}{i!}\left(x_{R}-X_t\right)^i\right)\D N_t\\
    &= f^{'}(X_t)\sigma_t\D W_t+\left( f^{'}(X_t) \mu_t+f^{''}(X_t) \frac{\sigma_t^2}{2} \right)\D t+\left( \sum\displaylimits_{i=1}^\infty\frac{f^{(i)}(X_t)}{i!}\left(x_{R}-X_t\right)^i+f(X_t)-f(X_t)\right)\D N_t\\
    &=f^{'}(X_t)\sigma_t\D W_t+\left( f^{'}(X_t)\mu_t+f^{''}(X_t) \frac{\sigma_t^2}{2} \right)\D t+\left( \sum\displaylimits_{i=0}^\infty\frac{f^{(i)}(X_t)}{i!}\left(x_{R}-X_t\right)^i -f(X_t)\right)\D N_t.
\end{align*}
Finally, using the Taylor formula we get the generalised It\^o formula \begin{equation*}
    \D Z_t=f^{'}(X_t)\sigma_t\D W_t+\left(f^{'}(X_t)\mu_t+f^{''}(X_t) \frac{\sigma_t^2}{2}  \right)\D t+\left( f(x_{R}) -f(X_t)\right)\D N_t.
\end{equation*}

For $\mu_t = 0$ and $\sigma_t =1$ we recover the It\^o lemma for the resetting process \eqref{eq:mysde}
 \begin{align}\label{lem:myito}
\D f(X_t)=\frac{1}{2}f^{''}(X_t)\D t+f^{'}(X_t)\D W_t+\left( f(x_{R})-f(X_t) \right)\D N_t.
\end{align}
\subsection{Moment generating function and Fourier transform}
Now, using the above It\^o lemma, let us calculate the MGF $M_t(s)=\mathbb{E}\left( e^{sX_t}\right)$ of the resetting process $X_t$ given by \eqref{eq:mysde}. Here $\E$ is the expected value. 
Applying \eqref{lem:myito} for the function $f_s(x)=e^{sx}$ and putting $Z_t=f_s(X_t)$ we arrive at the following SDE:
\begin{equation}\label{eq:sdemgf}
\D Z_t=\D f_s(X_t) = sZ_t\D W_t+\frac{1}{2}s^2Z_t\D t+\left( f_s(x_{R})-Z_t \right)\D N_t.
\end{equation}
Note that $\mathbb{E}\left(Z_t \right)$ is equal to the MGF of  $X_t$.
Taking the expected value of both sides of \eqref{eq:sdemgf} we get
\begin{align*}
    \E\left( \D Z_t \right)&=s{\E\left(  Z_t\D W_t \right)}+\frac{1}{2}s^2\D t\E( Z_t)+\E\left( e^{sx_{R}}-f(X_t) \right)\D N_t\\
    &=\frac{1}{2}s^2\D t\E (Z_t)+\sum\displaylimits_{i=0}^1\E\left(\left( e^{sx_{R}}-f(X_t) \right)\D N_t\Big{|}\D N_t=i\right)\PR \left( \D N_t = i \right)\\
    &=\frac{1}{2}s^2\D t\E( Z_t)+\E\left(e^{sx_{R}}-f(X_t)\Big{|}\D N_t=1\right)\PR \left( \D N_t = 1 \right)\\
    &=\frac{1}{2}s^2\D t\E (Z_t)+r\D te^{sx_{R}}-r\D t\E \left( f(X_t)\right).
\end{align*}
The sum in the second line of the above equalities is the consequence of the total probability formula.
Dividing both sides by $\D t$, interchanging the order of integration and differentiation and using the dominated convergence theorem (the process has clearly lighter tails than the ordinary diffusion, and as such the theorem can be applied), substituting $\E (Z_t)=M_t(s)$ and rearranging the equation, we get a simple ordinary differential equation
\begin{gather*}
    \partial_tM_t(s)+\left(r-\frac{1}{2}s^2\right)M_t(s)=re^{sx_{R}}.
\end{gather*}
This yields the solution
\begin{gather*}
    M_t(s)=\frac{re^{sx_{R}}}{ r-\frac{1}{2}s^2}+c(s)e^{-\left( r-\frac{1}{2}s^2 \right)t}.
\end{gather*}
Using conditions $M_0(s)=e^{sX_0}$ and $M_t(0)=1$ for all $t$ we get that the unique constant $c(s)$ equals
\begin{equation*}
    c(s)=\left( e^{sX_0}-\frac{r}{r-\frac{1}{2}s^2}e^{sx_{R}} \right).
\end{equation*}
\begin{figure}[htb]
    \includegraphics[scale=0.35]{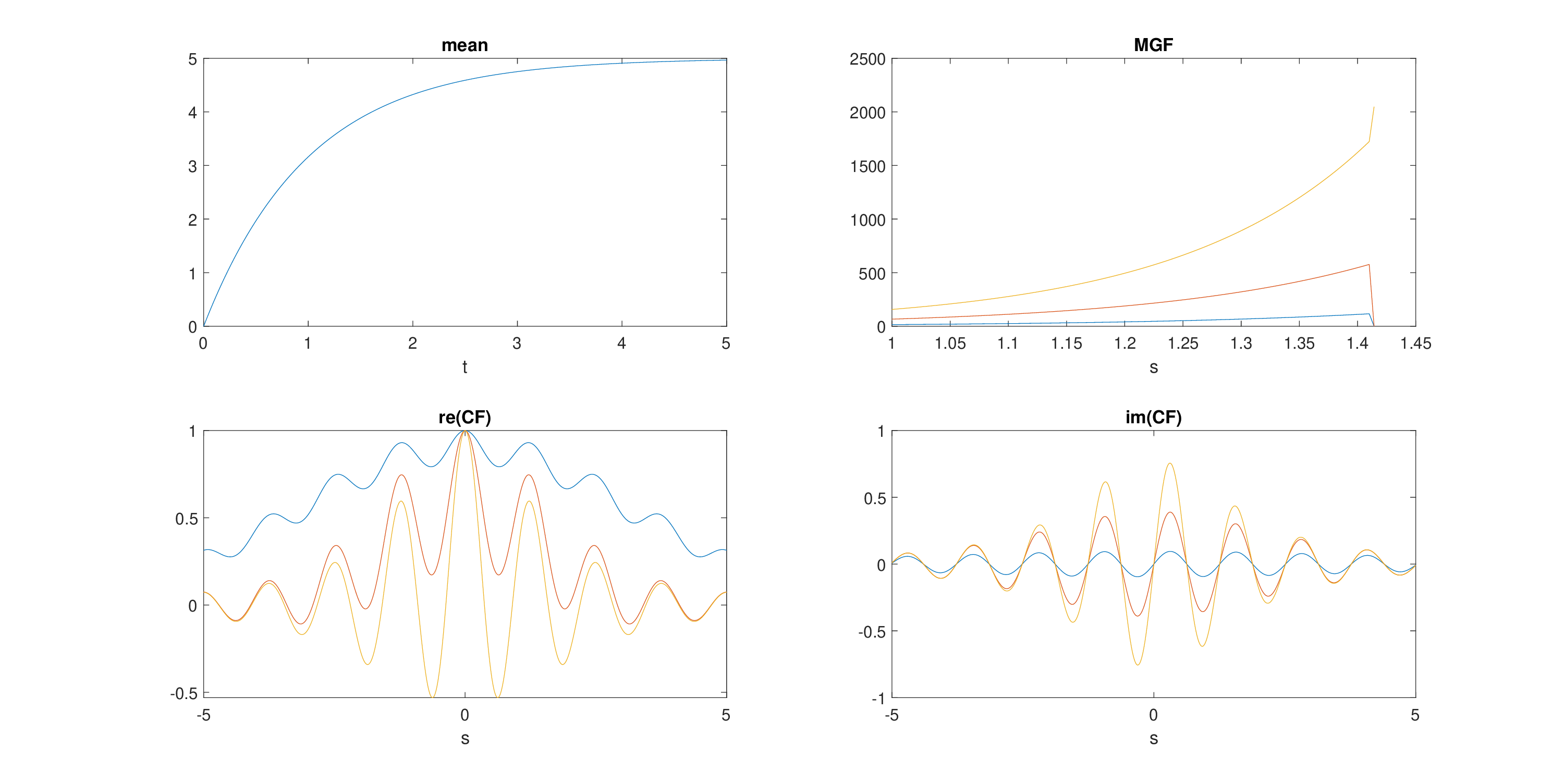}
    \caption{A plot of the derived mean of the resetting process (top left), MGF (top right), the real part of the Fourier transform (bottom left) and the imaginary part of the Fourier transform (bottom right) at a different time points (0.1, 0.5, 1.5), as a blue, red and yellow lines respectively. The process starts at 0 and resets with rate 1 to point 5. We can also see the unusual behavior of the MGF at the point $\sqrt{2}$, due to the divergence of the function.}
    \label{fig:analyticalresults}
\end{figure}
Finally our MGF has the form
\eqb
\label{MGF}
        M_t(s)=\frac{re^{sx_{R}}}{ r-\frac{1}{2}s^2}+\left( e^{sX_0}-\frac{re^{sx_{R}}}{r-\frac{1}{2}s^2} \right)e^{-\left( r-\frac{1}{2}s^2 \right)t}.
\eqe
We can also notice that in case of the lack of resetting ($r=0$) the MGF simplifies to the MGF of Brownian motion starting at $X_0$
\begin{equation*}
    \frac{0e^{sx_{R}}}{ 0-\frac{1}{2}s^2}+\left( e^{sX_0}-\frac{0}{0-\frac{1}{2}s^2}e^{sx_{R}} \right)e^{-\left( 0-\frac{1}{2}s^2 \right)t}=e^{sX_0+t\frac{1}{2}s^2}.
\end{equation*}
Additionally for large time $t\rightarrow\infty$  MGF of the process with resetting converges to
\begin{equation*}
    M(s)=\frac{e^{sx_{R}}}{1-\frac{1}{2r}s^2},
\end{equation*}
which is a MGF of the Laplace distribution with mean $x_{R}$ and scale parameter $(2r)^{-\frac{1}{2}}$. 

The Fourier transform of the process with resetting $ \varphi_t(s)=\mathbb{E}\left( e^{isX_t}\right) $ can be derived analogously or using the fact that $ \varphi_t(s)=M_t(is)$, and equals
\eqb
\label{FT}
        \varphi_t(s)=\frac{re^{isx_{R}}}{ r+\frac{1}{2}s^2}+\left( e^{isX_0}-\frac{re^{isx_{R}}}{r+\frac{1}{2}s^2} \right)e^{-\left( r+\frac{1}{2}s^2 \right)t}.
\eqe
Plots of $M_t(s)$ and $\varphi_t(s)$  are presented in Fig. \ref{fig:analyticalresults}. We can immediately notice that the MGF is well defined only for $|s|<\sqrt{2r}$. 
The moment generating function and Fourier transform are closely related via formula 
$\varphi_t(s)=M_t(is)$. However, the advantage of $\varphi_t(s)$ is that it is well defined for any $s\in\mathbb{R}$. The moment generating function does not have this property. On the other hand, the advantage of $M_t(s)$ is that it is real not complex. Both functions $M_t(s)$ (if well defined) and $\varphi_t(s)$ determine uniquely the underlying distribution.

\subsection{Explicit PDF}
\begin{figure}[hb]
    \includegraphics[scale=0.3]{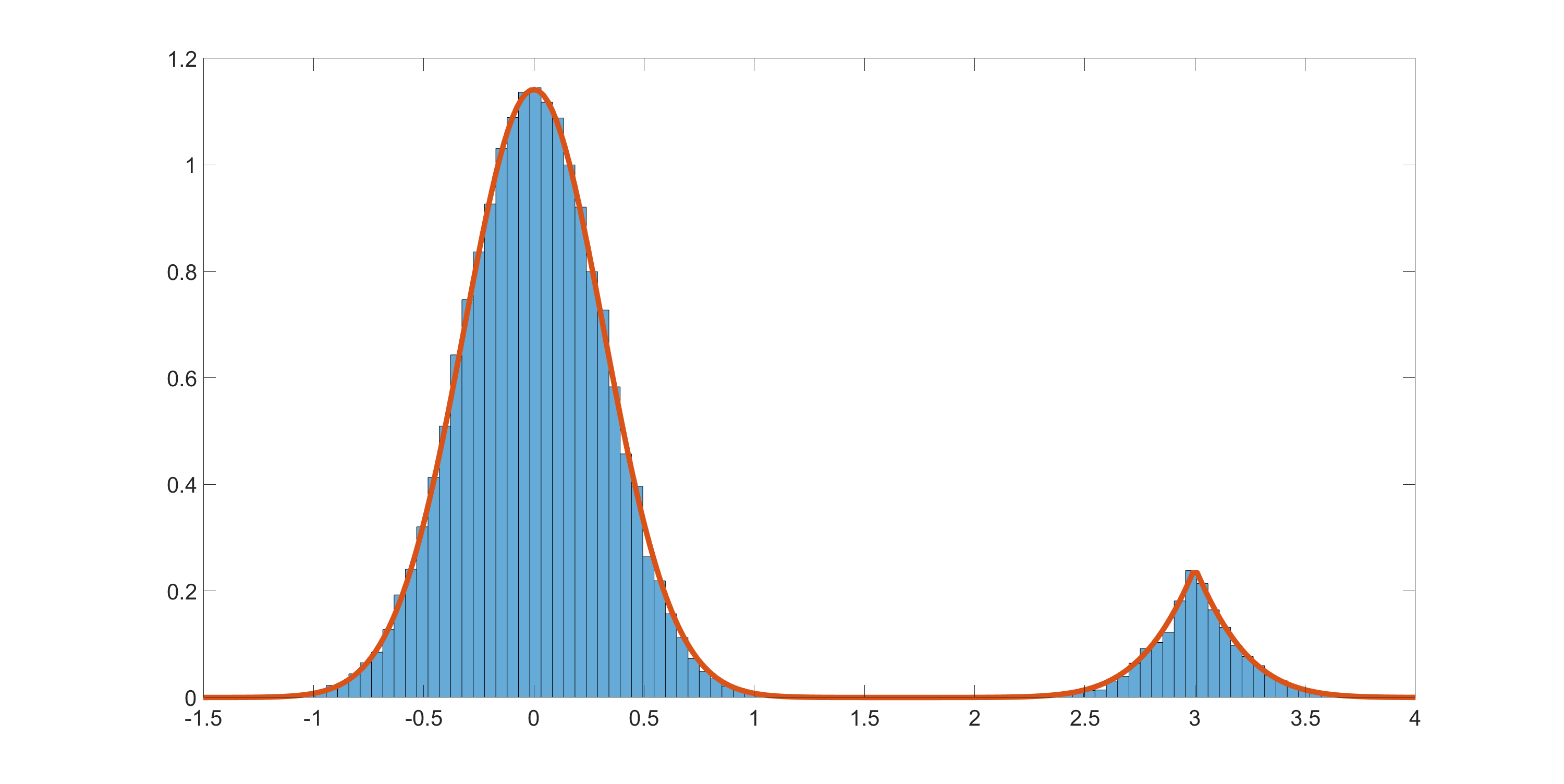}
    \caption{A histogram of $X_t$ obtained using Monte Carlo simulations compared with the derived analytical PDF of the process at $t=0.1$. We observe perfect agreement between both results. The process starts at 0 and resets to point 3 with rate 1.}
    \label{fig:analyticalpdf}
\end{figure}

Now, using the inverse Fourier transform we are able to invert $\varphi_t(s)$. This way we will arrive at the explicit formula for the PDF $p(x,t)$ of the resetting process $X_t$. An exemplary plot of $p(x,t)$ is presented in Fig. \ref{fig:analyticalpdf}.
Let us split the Fourier transform \eqref{FT} into three parts in the following way
\begin{equation*}
    \varphi_t(s)=\varphi_1+\varphi_2+\varphi_3,
\end{equation*}
where $\varphi_1=\frac{r}{ r+\frac{1}{2}s^2}e^{isx_{R}}$, $\varphi_2=e^{isX_0}e^{-\left( r+\frac{1}{2}s^2 \right)t}$, $\varphi_3=-\frac{r}{r+\frac{1}{2}s^2}e^{isx_{R}} e^{-\left( r+\frac{1}{2}s^2 \right)t}$. Using the well known properties of the Fourier transform, we can already see that the first term is just the Fourier transform of the Laplace distribution with the corresponding PDF \cite{Bill86}
\eqb
\label{Laplace_PDF}
 \F^{-1}\left\{ \varphi_1 \right\}=f_\mathcal{L}(x) = \sqrt{\frac{r}{2}} \exp\left( -\frac{|x-x_{R}|}{(2r)^{-1/2}} \right),
\eqe
see also \cite{reed2002gene,reed2003pareto} for the relationship between Laplace distribution and stopping times of random motions.

We will calculate the second term directly from the definition of the inverse transform
\begin{align*}
   \F^{-1}\left\{ \varphi_2 \right\}&=\frac{1}{2\pi}\int_{\mathbb{R}}\displaylimits e^{-isx}e^{isX_0-rt-\frac{t}{2}s^2}\D s=\frac{e^{-rt}}{2\pi}\int_{\mathbb{R}}\displaylimits e^{i(X_0-x)s-\frac{t}{2}s^2}\D s\\
    &=   \frac{e^{-rt-\frac{(x-X_0)^2}{2t}}}{\sqrt{2\pi t}}=
    e^{-rt}f_{W_{X_0}} (x),
\end{align*}
where $f_{W_{X_0}} (x)=\frac{e^{-\frac{(x-X_0)^2}{2t}}}{\sqrt{2\pi t}}$ is the PDF of standard Brownian motion starting at $X_0$.

Let us now write $\varphi_3$ in the following way
\begin{equation*}
    \varphi_3=-e^{-rt}\cdot\frac{r}{r+\frac{1}{2}s^2}e^{isx_{R}}\cdot e^{-\frac{t}{2}s^2}.
\end{equation*}
The first term of the product is independent of $s$, so it's treated as a scalar in the inverse Fourier transform. The second one is the Fourier transform of the Laplace distribution and the third is simply the Fourier transform of the normal distribution $\mathcal{N}(0,t)$. Having all this and using the well-known fact that  $\F\{f\}\F\{g\}=\F\{f*g\}$, where $f*g$ is the convolution of functions,
 we finally get that the explicit PDF $p(x,t)$ of the resetting process \eqref{eq:mysde} equals
\begin{equation}
\label{eq:mypdf}
    p(x,t)=f_\mathcal{L}(x)+e^{-rt}\left(f_{W_{X_0}}(x)-f_{\mathcal{N}(0,t)}* f_\mathcal{L}(x)\right).
\end{equation}
Here $f_{\mathcal{N}(0,t)}(x)$ is the PDF of $\mathcal{N}(0,t)$. Note that if $f$ and $g$ are two PDFs  
of some independent random variables, say V and W, then $f*g$ is the PDF of $V+W$. 
We underline that the convolution of Laplace and normal distributions can be calculated explicitly, see  \cite{geraci2017mixedeffects}.

\subsection{Moments}

Applying the above formula for PDF we are also able to calculate the mean of the resetting process
\begin{equation*}
    \E (X_t)=\int_{\mathbb{R}} xp(x,t)\D x=\int_{\mathbb{R}} xf_\mathcal{L} (x)\D x+ e^{-rt}\left(\int_{\mathbb{R}} x\left(f_{W_{X_0}}(x)-f_{\mathcal{N}(0,t)}*f_\mathcal{L}(x)\right)\D x\right) .
\end{equation*}
Using the fact that the convolution of PDFs is the PDF of the sum of the corresponding independent random variables, we get that
\begin{equation*}
    \E (X_t) = x_{R}+e^{-rt}\left(X_0- x_{R}\right).
\end{equation*}
Let us now derive the formula for the $n$-th moment of the resetting process. For simplicity of the notation we will focus on the case $x_{R}=0$. We have 
\begin{equation*}
    \E (X_t^n)=\int_{\mathbb{R}}\displaylimits x^np(x,t)\D x=\int_{\mathbb{R}}\displaylimits x^n f_\mathcal{L}(x)\D x + e^{-rt}\left( \int_{\mathbb{R}}\displaylimits x^nf_{W_{X_0}}(x)\D x-\int_{\mathbb{R}}\displaylimits x^nf_{\mathcal{N}(0,t)}*f_\mathcal{L}(x)\D x \right).
\end{equation*}
The first integral is the $n$-th moment of the Laplace random variable $L$ with mean $0$ and scale parameter $(2r)^{-\sfrac{1}{2}}$. Using the scaling property it can be rewritten as $L=(2r)^{-\sfrac{1}{2}}Y$, where $Y\sim\mathcal{L}(0,1)$. It is well known, that the $n$-th moment of $Y$ equals $n!$ for even $n$ and otherwise zero. Therefore
\begin{equation}
\label{moments_L}
    \E (L^n) =\E \left(((2r)^{-1/2}Y^n)\right)=(2r)^{-n/2}\E Y^n=(2r)^{-n/2}n!
\end{equation}
for even $n$ and zero otherwise.
The second integral being the $n$-th moment of the normal distribution $\mathcal{N}(X_0,{t})$ can be stated as follows \cite{winkelbauer2014moments}
\begin{align} 
\label{moments_N}
 \E (W^n_{X_0})= &=
 \begin{cases}
(\sqrt{2t})^{n}\frac{\Gamma((n+1)/2)}{\sqrt{\pi}} \Phi\left(-\frac{n}{2},\frac{1}{2};-\frac{X_0^2}{2t}\right) &\text{ for even } n \\
X_0(\sqrt{t})^{n-1}2^{(n+1)/2}\frac{\Gamma(n/2+1)}{\sqrt{\pi}} \Phi\left(\frac{1-n}{2},\frac{3}{2};-\frac{X_0^2}{2t}\right) &\text{ for odd } n.
 \end{cases}
 \end{align}
Here $\Phi(a,b;c)$ is the Kummer’s confluent hypergeometric function  \cite{winkelbauer2014moments}.

The case $X_0=0$ simplifies the results even more
\begin{equation}
    \E (W^n_0)=\E (W^n )=t^{\frac{n}{2}}(n-1)!!
\end{equation}
for even $n$ and zero otherwise.

The third integral is the most interesting. As we can see, it represents the $n$-th moment of the sum of normal random variable $W\sim\mathcal{N}(0,{t})$ and Laplace random variable $L\sim\mathcal{L}(0,(2r)^{-\sfrac{1}{2}})$. Note that $W$ and $L$ are independent. Let us calculate the $n$-th moment of their sum
\begin{equation*}
    \int\displaylimits_{\R}x^n f_{\mathcal{N}(0,t)}*f_\mathcal{L}\D x=\E ((W+L)^n).
\end{equation*}
Utilizing the binomial theorem we get the following result
\begin{align}
\nonumber\E ((W+L)^n) &=    \E \left( \sum\displaylimits_{k=0}^n \binom{n}{k}W^kL^{n-k}\right)=\sum\displaylimits_{k=0}^n \binom{n}{k}\left(\E( W^k)\right) \left(\E (L^{n-k})\right)\\
&=
\sum\displaylimits_{\substack{k=0,\\k\text{ even}}}^n \binom{n}{k}\left(\E (W^k)\right) \left(\E (L^{n-k})\right),
\label{moments_sum}
\end{align}
for even $n$ and zero otherwise.
This is due to the fact that that both $W$ and $L$ are independent and symmetric with mean $0$.

Gathering all the terms \eqref{moments_L},\eqref{moments_N} and \eqref{moments_sum} we get the final expression for the $n-th$ moment of $X_t$
\begin{equation*}
    \E( X_t^n) = \E( L^n)+e^{-rt}\E( W_{X_0}^n)- e^{-rt}\E\left((W+L)^n\right).
\end{equation*}

We note that another way of deriving moments of resetting process is by calculating derivatives of the MGF $M_t(s)$ in eq. \eqref{MGF}.

\subsection{Fokker-Planck equations}

Having the PDF \eqref{eq:mypdf} of the resetting process we can directly derive the corresponding Fokker-Planck equation. Let us calculate the necessary derivatives of \eqref{eq:mypdf}:
\begin{align*}
    \partial_t p(x,t)&= \partial_t\left( f_\mathcal{L}+e^{-rt}\left(f_{W_{X_0}}-f_{\mathcal{N}(0,t)}*f_\mathcal{L}\right)\right)=\partial_t\left( e^{-rt}\left(f_{W_{X_0}}-f_{\mathcal{N}(0,t)}*f_\mathcal{L}\right)\right)\\
    &=-re^{-rt}\left(f_{W_{X_0}}-f_{\mathcal{N}(0,t)}*f_\mathcal{L}\right)+e^{-rt}\partial_t\left(f_{W_{X_0}}-f_{\mathcal{N}(0,t)}*f_\mathcal{L}\right)\\
    &=-re^{-rt}\left(f_{W_{X_0}}-f_{\mathcal{N}(0,t)}*f_\mathcal{L}\right)+e^{-rt}\left( \partial_tf_{W_{X_0}}-\left( \partial_t f_{\mathcal{N}(0,t)} \right)*f_\mathcal{L} \right)\\
    &=-rp(x,t)+rf_\mathcal{L}+e^{-rt}\left( \partial_tf_{W_{X_0}}-\left( \partial_t f_{\mathcal{N}(0,t)} \right)*f_\mathcal{L} \right)
\end{align*}
The second spatial derivative equals to
\begin{align*}
    \partial_{xx}p(x,t)&=\partial_{xx}\left( f_\mathcal{L}+e^{-rt}\left(f_{W_{X_0}}-f_W*f_\mathcal{L}\right)\right)=\partial_{xx}f_\mathcal{L}+e^{-rt}\left( \partial_{xx}f_{W_{X_0}}-\left( \partial_{xx}f_W \right)*f_\mathcal{L} \right).
\end{align*}
We can now note that the PDF of the Wiener process fulfills the standard diffusion equation 
\begin{equation*}
    \partial_t f_{W_{X_0}}=\frac{1}{2}\partial_{xx} f_{W_{X_0}},
\end{equation*}
so the time derivative can be rewritten as
\begin{align*}
    \partial_tp(x,t)&=-rp(x,t)+rf_\mathcal{L}+\frac{1}{2}\left(\partial_{xx}f_\mathcal{L} +e^{-rt}\left( \partial_{xx}f_{W_{X_0}}-\left( \partial_{xx} f_{\mathcal{N}(0,t)} \right)*f_\mathcal{L} \right)-\partial_{xx}f_\mathcal{L}\right)\\
    &=\frac{1}{2}\partial_{xx}p(x,t)+\left(r-\frac{1}{2}\partial_{xx} \right)f_\mathcal{L}-rp(x,t)
\end{align*}
The derivative (in the weak sense) of the density of Laplace distribution is
\begin{align*}
    \partial_xf_\mathcal{L}=\partial_x\sqrt{\frac{r}{2}}e^{-\sqrt{2r}|x-x_{R}|}=-re^{-\sqrt{2r}|x-x_{R}|}\partial_x\left( |x-x_{R}| \right)=-re^{-\sqrt{2r}|x-x_{R}|}\sgn(x-x_{R}),
\end{align*}
while the second spatial derivative is
\begin{align*}
    \partial_{xx}f_\mathcal{L}&=\sqrt{2r^3}e^{-\sqrt{2r}|x-x_{R}|}\sgn^2(x-x_{R})-re^{-\sqrt{2r}|x-x_{R}|}2\delta(x-x_{R})\\
    &=e^{-\sqrt{2r}|x-x_{R}|}\left(2r\delta(x-x_{R}) -\sqrt{2r^3}\sgn^2(x-x_{R}) \right)\stackrel{\text{a.e.}}{=}f_\mathcal{L}\left( 2\sqrt{2r}\delta(x-x_{R})-2r\right).
\end{align*}
Plugging this result to the expression for time derivative we arrive at the final Fokker-Planck equation for the resetting process
\begin{equation*}
    \partial_tp(x,t)=\frac{1}{2}\partial^2_xp(x,t)+\sqrt{2r}\delta(x-x_{R})f_\mathcal{L}-rp(x,t).
\end{equation*}
The above equation is different from \cite{Evans2020}, see \eqref{eq:evanspdf}, however we can check the equivalence of both equations by looking at them in the frequency domain via Fourier transform. We only need to look at the component $\sqrt{2r}\delta(x-x_{R})f_\mathcal{L}$. We have
\begin{align*}
    \F \left\{ \sqrt{2r}\delta(x-x_{R})f_\mathcal{L} \right\}&=\sqrt{2r}\int\displaylimits_{\R}\sqrt{\frac{r}{2}}e^{-\sqrt{2r}|x-x_{R}|}\delta(x-x_{R})e^{isx}\D x\\
    &=re^{-\sqrt{2r}|x_{R}-x_{R}|}e^{isx_{R}}=re^{isx_{R}}=\F\left\{ r\delta(x-x_{R}) \right\}.
\end{align*}
Both equations are equal in the frequency domain, implying that their solutions coincide.
One should also underline that the difference between both equations is only in the Dirac delta term
$\delta(x-x_{R})$, which affect the solutions only on the set of points with zero Lebesgue measure.
Since the random variable $X_t$ is continuous, modification of its PDF on the set of points with zero Lebesgue measure does not change the distribution.

Let us note that the above Fokker-Planck equation can be used to derive the stationary PDF of the resetting process just by putting $\partial_tp(x,t)=0$ and calculating the corresponding $p$.
\subsection{Infinitesimal generator}

The infinitesimal generator is a key operator in the theory of Markov processes \cite{ethier2009markov}.
It contains great deal of information about the Markov process. In particular it can be used to find the corresponding evolution equations. Applying the general formula for jump diffusion processes \cite{OksendalSulem}, we get that the infinitesimal generator $\mathcal{A}$ of the resetting process \eqref{eq:mysde} equals:

\begin{equation}\label{eq:gen}
    \mathcal{A}g(x)=\frac{1}{2}\partial_{xx}g(x)+r\left( g(x_{R})-g(x)\right),
\end{equation}
for appropriately smooth function $g$. The first component $\frac{1}{2}\partial_{xx}g(x)$ is the diffusive part of the generator. The second component (jump part) of the above expression can be identified as a generator of the compound Poisson process with jump sizes $x_{R}-X_t$, which corresponds to the jumps of the process to the resetting point $x_{R}$ with intensity $r$.

Let us now derive the adjoint generator of our resetting process. It follows from the calculations below
\begin{align*}
    \int f(z)\mathcal{O}g(z)\D z&=\int g(x_{R})f(z)\D z=\int g(u)\delta(u-x_{R})\D u\int f(z)\D z\\
    &=\iint g(u)\delta(u-x_{R})f(z)\D z \D u=\int g(u)\delta(u-x_{R})\int f(z)\D z \D u\\
    &=\int g(z)\delta(z-x_{R})\int f(u)\D u \D z=\int g(z)\mathcal{O}^*f(z)\D z,
\end{align*}
that the adjoint generator $\mathcal{A}^*$ equals
\begin{equation}
\label{adjoint}
    \mathcal{A}^*g(x)=\frac{1}{2}\partial_{xx}g(x)+r\left( \delta(x-x_{R})\int\displaylimits_{\mathbb{R}} g(z)\D  z-g(x)\right).
\end{equation}
From the general theory of Markov processes we know that the PDF of a Markov process with adjoint generator $\mathcal{A}^*$ satisfies the following Fokker-Planck equation \cite{pavliotis2014stochastic}
\begin{equation*}
\partial_t p=\mathcal{A}^* p.
\end{equation*}
Therefore, applying \eqref{adjoint} we obtain the following Fokker-Planck formula for the resetting process
\begin{equation*}
     \partial_tp(x,t)=\frac{1}{2}\partial_{xx}p(x,t)+r\left( \delta(x-x_{R})\int\displaylimits_{\mathbb{R}} p(z,t)\D  z-p(x,t)\right).
\end{equation*}
Since $\int\displaylimits_{\mathbb{R}} p(z,t)\D  z =1$ we obtain the same Fokker-Planck equation  as the one derived in \cite{Evans2020} for resetting processes, cf. equation \eqref{eq:evanspdf}.

\subsection{Nonhomogeneous in time Poissonian resetting}

Non-homogeneous Poisson process (NPP) is a natural generalization of the standard Poisson process to the case of time-dependent intensity.  It is defined as a counting process with the following properties \cite{Ross1996}:
\begin{itemize}
    \item $\tilde{N}_0=0$,
    \item $\tilde{N}_t$ has independent increments,
    \item the increments $\tilde{N}_{t+h}-\tilde{N}_t$ are Poisson distributed with mean $\int_t^{t+h} r(s)ds$.

\end{itemize}
Here, the non-negative function $r(t)$ is called the intensity function of $\tilde{N}_t$. For higher values of $r(t)$ we observe more jumps of $\tilde{N}_t$, whereas small $r(t)$ gives less jumps on the average. For $r(t)=r=const$ we recover the standard Poisson process.

Now, we introduce a modified SDE defining our resetting process
\begin{equation}\label{eq:nonhomres}
    \D X_t=\D W_t+ (x_R-X_t)\D \tilde{N}_t.
\end{equation}
We put $X_0=x_{R}=0$ for simplicity. $\tilde{N}_t$ here is the NPP with intensity function $r(t)$ and $W(t)$ is the standard Brownian motion independent of $\tilde{N}_t$. 

In our further analysis we will focus on the power-law intensity function and the impact it has on both the mean square displacement (MSD) and the overall distribution. We will thus assume that the intensity function has the form $r(t)=r\cdot(t+1)^p$. The previously studied homogeneous resetting is now a special case with $p=0$, which makes it an interesting starting consideration point. With $p=0$ the process displayed a nondegenerate stationary distribution, so any increase in resetting frequency should have a profound impact on the probability law, most likely creating a degenerate distribution. On the other hand the decrease of the resetting intensity intuitively should recover the classical diffusion at infinity. This will be verified in detail. We will show that the diffusive, subdifussive, stationary and deterministic behavior can be observed for long times depending on the parameter $p$. 

Now let's define $R(t)=\int_0^tr(\omega)\D \omega$ as the mean number of resets up to point $t$. Using exactly the same methods and calculations as in the homogeneous case, based on It\^o lemma, we conclude that the Fourier transform of the process must be a solution of the following differential equation
\begin{equation*}
    \partial_t\varphi_t(s)+(r(t)+\frac{1}{2}s^2)\varphi_t(s)=r(t).
\end{equation*}
Integration yields the Fourier transform
\begin{equation}
\label{FT_NPP}
    \varphi_t(s)=e^{-R(t)-\frac{1}{2}ts^2}\left( 1+\int_0^tr(w)e^{R(w)+\frac{1}{2}ws^2}\D w \right)
\end{equation}
which after simple inspection can be inverted to obtain the density of the process given in \eqref{eq:nonhomres}
\begin{equation}\label{eq:nonhomresf}
    f_{X_t}(x,t)=e^{-R(t)}f_{\mathcal{N}(0,t)}(x)+e^{-R(t)}\int\displaylimits_0^tr(\omega)e^{R(\omega)}f_{\mathcal{N}(0,t-\omega)}(x)\D \omega.
\end{equation}
We can easily extract the MSD simply by multiplying by $x^2$ and integrating. We get
\begin{equation*}
    MSD(X_t)=e^{-R(t)}MSD(W_t)+e^{-R(t)}\int\displaylimits_0^tr(\omega)e^{R(\omega)}MSD(W_{t-\omega})\D \omega,
\end{equation*}
which after simple calculations turns into
\begin{equation}
\label{MSD_NPP}
    MSD(X_t)=e^{-R(t)}\int\displaylimits_0^te^{R(\omega)}\D \omega.
\end{equation}
Now, let us assume that the intensity function has the form $r(t)=r\cdot(t+1)^p$, where $p$ can be any real number and $r$ is a positive constant. Then we get $R(t)=\int_0^tr(\omega)\D \omega=\frac{r}{p+1}\left(\left(t+1)^{p+1}-1\right)\right)$ if $p\neq -1$ and $R(t)=\ln{(t+1)^r}$ if $p=-1$. 

We already know the asymptotic behavior of $X_t$ for $p=0$ from previous sections. Let us now check the  case $p=-1$. We have
\begin{equation*}
    MSD(X_t)=e^{-\ln{(t+1)^r}}\int\displaylimits_0^te^{\ln{(\omega+1)^r}}\D \omega=\frac{1}{(t+1)^r}\int\displaylimits_0^t(\omega+1)^r\D \omega
\end{equation*}
The MSD is now equal to $\frac{t+1}{r+1}-\frac{1}{(t+1)^r(r+1)}$, which for large times gives us linear diffusive scaling. 

Next, if $p\in(-1,0)$ then $r(t)$ is decreasing, but $R(t)$ is increasing, giving us the most interesting case to investigate. We obtain
\begin{equation*}
    MSD(X_t)=e^{-\frac{r}{p+1}\left(\left(t+1)^{p+1}-1\right)\right)}\int\displaylimits_0^te^{\frac{r}{p+1}\left(\left(\omega+1)^{p+1}-1\right)\right)}\D \omega=e^{-\frac{r}{p+1}(t+1)^{p+1}}\int\displaylimits_0^te^{\frac{r}{p+1}(\omega+1)^{p+1}}\D \omega.
\end{equation*}
Applying the de l'Hospital rule we get for large times
\begin{equation*}
    MSD(X_t)=\frac{\int\displaylimits_0^te^{\frac{r}{p+1}(\omega+1)^{p+1}}\D \omega}{e^{\frac{r}{p+1}(t+1)^{p+1}}}\approx\frac{e^{\frac{r}{p+1}(t+1)^{p+1}}}{r(t+1)^pe^{\frac{r}{p+1}(t+1)^{p+1}}}\approx\frac{1}{r}(t+1)^{-p}.
\end{equation*}
The process now displays subdiffusive behavior with exponent $-p\in(0,1)$.

Consequently, for $p<-1$ we obtain for large times
\begin{equation*}
    MSD(X_t)=e^{-\frac{r}{p+1}(t+1)^{p+1}}\int\displaylimits_0^te^{\frac{r}{p+1}(\omega+1)^{p+1}}\D \omega\approx \int\displaylimits_0^te^{\frac{r}{p+1}(\omega+1)^{p+1}}\D \omega
\end{equation*}
Again, applying the de l'Hospital rule we get 
\begin{equation*}
    MSD(X_t)\approx t
\end{equation*}
for large times. As expected, the process displays normal diffusive scaling of the MSD. 

Examining the case $p>0$ in analogous way, we find that  $MSD(X_t)\approx \frac{1}{r}(t+1)^{-p}$, which converges to 0 for large times. 

Let us now focus on the asymptotic distributions.  Again, we expect different results for different $p$. The case $p=0$ corresponds to the homogeneous Poisson process and was analyzed previously, giving rise to the Laplace stationary distribution.

For $p>0$ we get from \eqref{FT_NPP} and from de l'Hospital rule that the Fourier transform satisfies  
\[
\varphi_t(s) \rightarrow 1
\] 
for large times. This implies that $X_t\stackrel{d}{\rightarrow} x_R = 0$ as $t\rightarrow \infty$. Here $\stackrel{d}{\rightarrow}$ means convergence in distribution. This result agrees with our intuition. For $p>0$ the number of resetting times increases with time, therefore $X_t$ keeps coming back to $x_R$ more and more often.

For $p\in(-1,0)$ we get from \eqref{FT_NPP} and from de l'Hospital rule that the Fourier transform satisfies  
\[
\varphi_t(s) \approx \frac{1}{1+\frac{s^2}{2r(t+1)^p}}
\] 
for large times.
This means that the distribution of $X_t$ is asymptotically Laplace with subdiffusive $MSD(X_t)\approx t^{-p}$ for large times.

Ananlogous result is obtained for $p=-1$:
\begin{equation*}
    \varphi_{t}(s)\approx \frac{1}{1+\left(\frac{t+1}{2r}\right)s^2}    
\end{equation*}
for large times, with linear in time asymptotic MSD.

\begin{figure}[t!]
    \centering
    \includegraphics[scale=0.3]{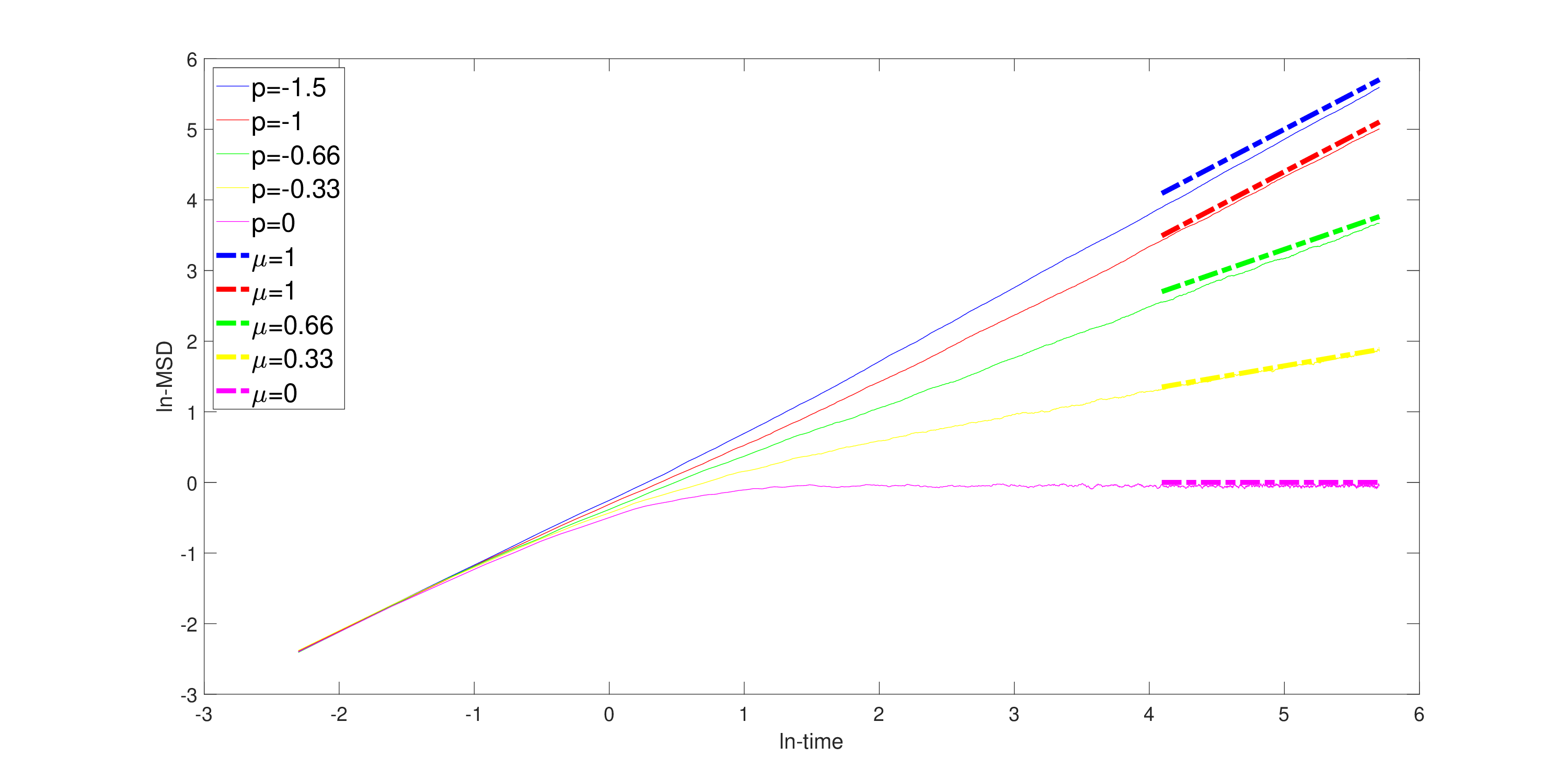}
    \caption{A Monte-Carlo simulation results (solid lines) compared with the theoretical asymptotics  (dotted lines) of the MSD of Brownian motion under nonhomogeneous Poisson resetting with  $r(t)=r\cdot(t+1)^p$, calculated for different $p\leq 0$. The parameter $\mu$ here is the corresponding power-law exponent of the theoretical asymptotic MSD. 
The time horizon of the simulations was of order $10^2$, while the number of samples was equal to $10^4$. We can observe the convergence of the anomalous diffusion exponents to the  theoretical ones.}
    \label{fig:diffusionexp}
\end{figure}
\begin{figure}[t!]
    \centering
    \includegraphics[scale=0.3]{ 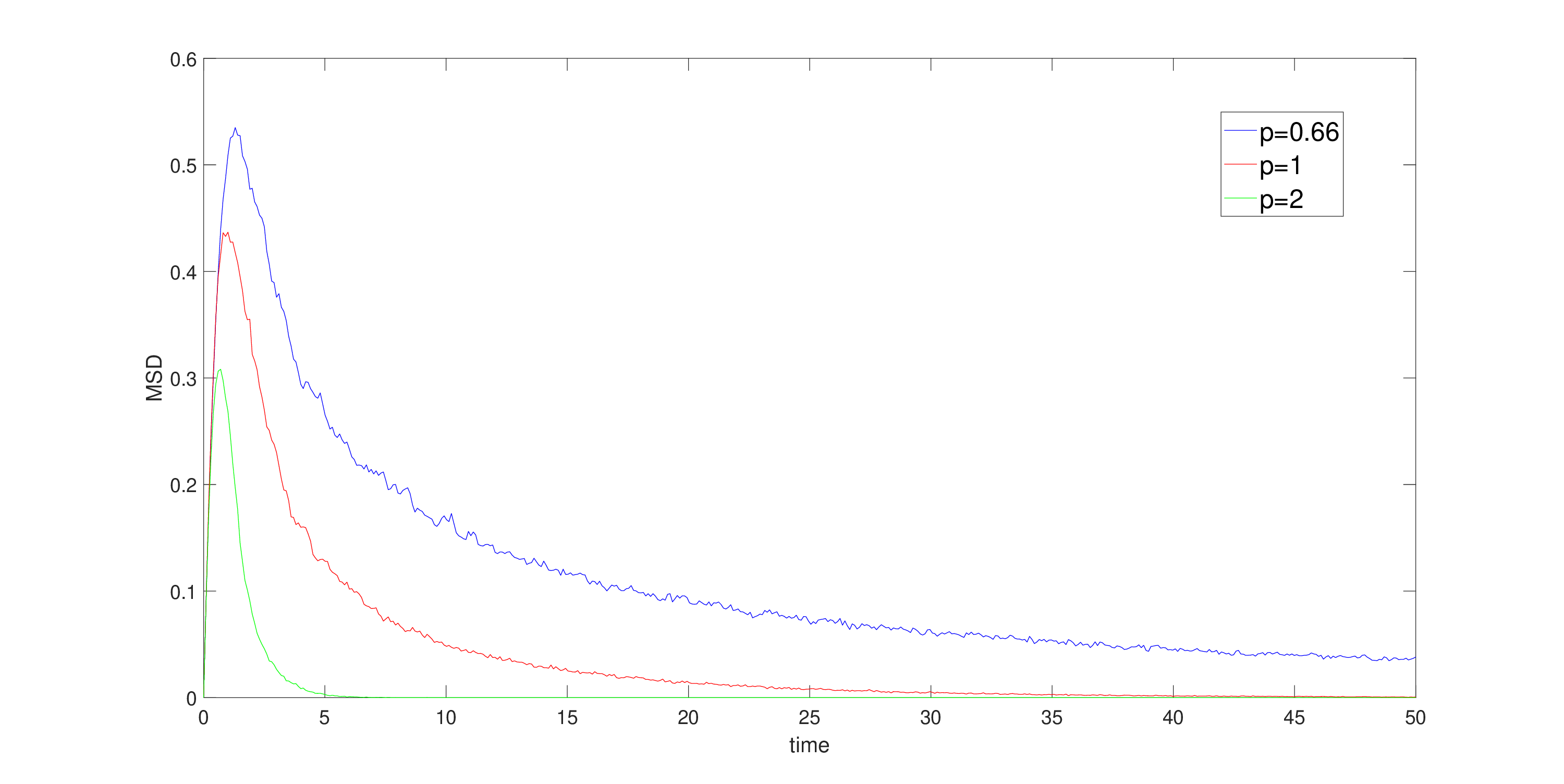}
    \caption{A Monte-Carlo simulation (solid lines) of the MSD of the process with nonhomogeneous Poisson resetting with intensity function  $r(t)=r\cdot(t+1)^p$, calculated for different $p>0$. 
Asterixes are the corresponding theoretical values of the MSD obtained via numerical evaluation of \eqref{MSD_NPP}.
We can see the convergence of the MSD to 0, indicating the convergence to the degenerate distribution concentrated in $x_R=0$. The convergence rate increases with the increase of the $p$ parameter.}
    \label{fig:diffusionexp0}
\end{figure}

Looking at $p<-1$ we get that for large times the corresponding characteristic function satisfies 
\[
\varphi_{t}(s)\approx e^{-R(t)}e^{-ts^2 /2}+ \frac{r(t)(t+1)}{p}\frac{1}{1+\frac{t+1}{2p}s^2} , 
\]
which is a combination of Gaussian and Laplace distributions. The corresponding MSD is asymptotically linear in time. Note that for $p\rightarrow -\infty$, i.e. when the resetting events vanish, only the diffusive Brownian part is left in the above formula. 

Summarizing the above findings for large times:
\begin{itemize}
    \item for $p>0$ we observe degenerate distribution concentrated in $x_{R}$ with $MSD(X_t)\approx const$;
    \item for $p=0$ a stationary Laplace distribution with $MSD(X_t)\approx const$ is observed;
    \item for $p\in(-1,0)$ we get a non-stationary Laplace distribution with subdiffusive $MSD(X_t)\approx t^{-p}$;
 \item for $p=-1$ we get a non-stationary Laplace distribution with diffusive $MSD(X_t)\approx t$;
    \item for $p<-1$ we observe a combination of normal and Laplace distributions with diffusive $MSD(X_t)\approx t$. As  $p\rightarrow -\infty$, the standard Brownian diffusion is recovered.
\end{itemize} 

In Figs. \ref{fig:diffusionexp}--\ref{fig:diffusionexp0} we can observe the approximated MSD with a fitted (anomalous) diffusion exponent. We observe perfect agreement between simulation and theory.

\section{Conclusions}
Summarizing, we have introduce a general stochastic representation for processes with resetting. It allows to analyze such processes using stochastic differential equations, It\^{o} lemma and related tools. We have shown robustness of our approach for the case of homogeneous, and nonhomogeneous Poissonian resetting.
We would like to underline that the derived in this paper results build a general framework for the analysis of stochastic processes with intermittent random resetting. It allows to analyze processes with resetting both, analytically and using Monte Carlo simulation methods. 
Moreover, the presented approach can be easily generalized to the case of non-Markovian resetting times as well as other, not necessarily Brownian, driving processes (L\'evy flights, L\'evy processes, fractional processes, etc.).  

\begin{acknowledgments}
This research was partially supported by NCN Sonata Bis-9 grant nr 2019/34/E/ST1/00360.
\end{acknowledgments}
\bibliography{bibliografia} 
\end{document}